\documentclass[12pt]{article}
\usepackage{latexsym}
\usepackage{amssymb}
\usepackage{amsmath}
\usepackage{amsfonts}
\usepackage{amsthm}
\usepackage{picinpar}

\newtheorem{thm}{Theorem}
\newtheorem{lem}[thm]{Lemma}

\def\startproof{\bigskip \noindent {\em Proof }}
\def\endproof{\hfill $\vcenter{\hrule height .3mm
		\hbox {\vrule width .3mm height 2mm \kern 2mm
			\vrule width .3mm} \hrule height .3mm}$ \bigskip}

\begin{document}

\title{The functional equation for $\zeta$}
\author{Keith Ball}

%%%%%Acknowledgements

\maketitle

This note contains a short proof of the functional equation for the zeta function. The argument is computationally similar to the Hankel contour arguments found in many places but is made simpler to visualise by the symmetries of the sine. In this form the argument seems to be new or at least not widely known. We begin with a simple lemma which uses integration by parts to introduce the sine (or hyperbolic sine).

\begin{lem} For $s>1$
\[ \Gamma(s) \zeta(s)= \frac{1}{4 s} \int_0^{\infty} \frac{t^s}{\sinh^2 (t/2)} \, dt. \]
\end{lem}

\startproof 
\begin{eqnarray*}
\Gamma(s) \zeta(s) & = & \sum_1^{\infty} \frac{1}{n^s} \int_0^{\infty} e^{-t} t^{s-1} \, dt = \sum_1^{\infty} \int_0^{\infty} e^{-n t} t^{s-1} \, dt \\
& = & \int_0^{\infty}  \frac{t^{s-1}}{e^t-1} \, dt = \frac{1}{s} \int_0^{\infty}  \frac{e^t \, t^{s}}{(e^t-1)^2} \, dt \\
& = & \frac{1}{4 s} \int_0^{\infty} \frac{t^s}{\sinh^2 (t/2)} \, dt.
\end{eqnarray*} \endproof

We now move on to the theorem itself.
\begin{thm} The function $s \mapsto (s-1)\zeta(s)$ defined on $(1,\infty)$ has an analytic continuation to the entire complex plane and for all $s$ we have 
\[  \zeta(s)=2 (2 \pi)^{s-1} \sin (\pi  s/2) \Gamma(1-s) \zeta(1-s) \]
with the usual convention regarding the removable singularities.
\end{thm}

\startproof Consider the integral
\[ \frac{1}{2 \pi i} \int_C \frac{\pi^2 z^{1-s}}{\sin^2 \pi z} \, dz \]
where the contour $C$ is a vertical line between 0 and 1, say $(\Re z=1/2)$, traversed upwards. In this and the rest of the argument we take powers $z^{-s}$ only for complex numbers satisfying $\Re z \geq 0$ so there is no ambiguity. Because of the rapid growth of $\sin \pi z$ as $z$ moves away from the real axis, this function is easily seen to be an entire function of $s$.

If $s>1$ then the factor $(x+i y)^{1-s} \rightarrow 0$ uniformly in $y$ as $x \rightarrow \infty$ while the factor $\sin^2 \pi (x+i y)$ is periodic in $x$, so we can compute the integral as the negative of the sum of the residues of 
\[ \frac{\pi^2 z^{1-s}}{\sin^2 \pi z} \]
at the positive integers. The residue at $n$ is $(1-s)/n^s$ so we find that the integral is $(s-1)\zeta(s)$.

On the other hand, if $s<0$ the integral makes sense at $z=0$ so we may shift the contour onto the imaginary axis. By splitting into the upper and lower halves of the axis we get that the integral is 
\begin{eqnarray*}
 \frac{(e^{i \pi/2})^{2-s}-(e^{-i \pi/2})^{2-s} }{2 \pi i} \int_0^{\infty} \frac{\pi^2 y^{1-s}}{\sin^2 \pi i y} \, dy &=& -\pi \sin (\pi  s/2) \int_0^{\infty} \frac{y^{1-s}}{\sinh^2 \pi y} \, dy \\
 & = & \frac{-\pi \sin (\pi  s/2)}{(2 \pi)^{2-s}}  \int_0^{\infty} \frac{t^{1-s}}{\sinh^2 (t/2)} \, dt \\
 & = & \frac{-\pi \sin (\pi  s/2)}{(2 \pi)^{2-s}}  4(1-s) \Gamma(1-s) \zeta(1-s) \\
 & = & 2 (2 \pi)^{s-1} \sin (\pi  s/2) (s-1) \Gamma(1-s) \zeta(1-s) 
 \end{eqnarray*}
 where the last but one identity uses the lemma.
\endproof

\vspace{1in}

	\noindent Keith Ball \\
	Mathematics Institute \\
	University of Warwick \\
	Coventry CV4 7AL
	
	\noindent k.m.ball@warwick.ac.uk

\end{document}